\pdfoutput=1
\documentclass[12pt,twoside,reqno]{amsart}
\usepackage{amsfonts,amssymb,amsmath,amsthm,amsxtra,mathabx}
\usepackage{mathtools}
\usepackage{bm}
\usepackage{xparse}
\usepackage{mathtools}
\usepackage{dsfont}
\usepackage{bbm}
\usepackage{mathrsfs}
\usepackage{tikz-cd}
\usepackage{enumitem}
\usepackage{comment}
\usepackage{lmodern}
\usepackage[top=1.5in, bottom=1.5in, left=1.23in, right=1.23in]{geometry}
\usepackage[all]{xy}
\usepackage{amsmath}
\usepackage{amssymb}
\usepackage{amsthm}
\usepackage{amscd}

\DeclareMathAlphabet{\mathbfit}{OT1}{cmss}{bx}{it}
\DeclareMathAlphabet\mathbfcal{OMS}{cmsy}{b}{n}

\usepackage{scalerel,stackengine}
\stackMath
\newcommand\reallywidehat[1]{%
\savestack{\tmpbox}{\stretchto{%
  \scaleto{%
    \scalerel*[\widthof{\ensuremath{#1}}]{\kern-.7pt\bigwedge\kern-.7pt}%
    {\rule[-\textheight/2]{1.5ex}{\textheight}}
  }{\textheight}%
}{0.9ex}}%
\stackon[2.25pt]{#1}{\tmpbox}%
}

\DeclareMathAlphabet{\mathbfit}{OT1}{cmss}{bx}{it}
\DeclareMathAlphabet\mathbfcal{OMS}{cmsy}{b}{n}

\numberwithin{equation}{section}

\newtheorem{theorem}[equation]{Theorem}

\newtheorem{lemma}[equation]{Lemma}

\theoremstyle{definition}
\newtheorem{remark}[equation]{Remark}

\newtheorem{definition}[equation]{Definition}

\usepackage{hyperref} 
\hypersetup{
    colorlinks=true,       
    linkcolor=blue,          
    citecolor=magenta,        
    filecolor=magenta,      
    urlcolor=cyan           
}

\begin{document}
\title[A Variation norm Carleson theorem in higher dimensions]{A Variation norm Carleson theorem in higher dimensions}
\author[Himali Dabhi]{Himali Dabhi}
\address[Himali Dabhi]{School of mathematics, the university of edinburgh and the maxwell institute for the mathematical sciences, 
James Clerk Maxwell Building,
The King's Buildings,
Peter Guthrie Tait Road,
City Edinburgh,
EH9 3FD}
\email{h.p.dabhi@sms.ed.ac.uk}

\newcommand{\bk}[1]{{\color{red}{#1}}}

\setcounter{tocdepth}{1}
\renewcommand{\thefootnote}{}
\footnote{\emph{To appear in Colloquium Mathematicum}}
\renewcommand{\thefootnote}{\arabic{footnote}}
\setcounter{footnote}{0}

\date{}
\begin{abstract}
In 1971, C. Fefferman established a higher dimensional extension of the celebrated Carleson--Hunt theorem which gives pointwise almost everywhere convergence of the partial Fourier sums of functions in $L^p(\mathbb T), 1 < p < \infty.$ More precisely, Fefferman proved a maximal function bound for polygonal Fourier partial sums of functions in  $L^p(\mathbb T^d), p>1.$ In this note, we extend Fefferman's maximal function bound to strong $r$-variation norm bounds whenever $r>2$ as well as uniform $2$-oscillation bounds. Furthermore, for functions in $L^2(\mathbb T^d)$, we establish $r$-variational and $2$-oscillation bounds for partial Fourier sums over nested rectangles whenever $r>2$.
\end{abstract}

\maketitle

\section{Introduction}
Let $f$ be an integrable function on the torus $\mathbb T$ and let $\{\widehat{f}(n)\}_{n \in \mathbb Z}$ be its Fourier coefficients where $\widehat{f}(n)=\int_{\mathbb T} f(x)e^{-2 \pi i nx}dx$. We define the partial sums of the Fourier series of $f$ by
\begin{equation}
\begin{aligned}
    S_{N}f(x) : =\sum_{|n| \leq N} \widehat{f}(n)e^{2 \pi i n x}, \,\,\, N \in \mathbb N.
    \end{aligned}
\end{equation}
Carleson \cite{C} showed that for a square integrable function $f$, its Fourier series converges to $f$ almost everywhere on the torus $\mathbb T$. Hunt \cite{H} extended the result to all $L^p$ functions for $1 <p<\infty$ and proved the maximal function bound 
\begin{equation}\label{chthm}
\begin{aligned}
    \|\sup_{N \in \mathbb N}|S_{N}f|\|_{L^p(\mathbb T)} \leq C_p \|f\|_{L^p(\mathbb T)}.
    \end{aligned}
\end{equation}
In the last 30 or so years, there has been much interest in the ergodic and harmonic analysis communities to extend maximal function bounds (such as (\ref{chthm})) to variational and/or oscillation bounds; see for example, \cite{JSW} and \cite{MSW}. Such bounds immediately imply pointwise almost everywhere convergence and also give quantitative information on the rate of convergence. On the other hand, the passage from a maximal function bound to pointwise convergence requires knowledge of convergence for a dense family of functions. Thus variational and oscillation bounds are useful in settings where convergence for a dense family of functions is not clear. A general probability space is an example of such a setting which is also the setting for many results in ergodic theory.

For $1 \leq r < \infty$, the $r$-variation of a family $\{F_t:t \in \mathbb I\}$ of complex-valued measurable functions on a measure space $X$, where $\mathbb I \subseteq \mathbb R$, is defined by
\begin{equation}\label{rvardef}
\begin{aligned}
    V^r(F_t(x):t \in \mathbb I):=\sup_{t_1 \leq t_2 \leq \ldots \leq t_L} \Bigl( \sum_{l=1}^{L-1}|F_{t_{l+1}}(x)-F_{t_l}(x)|^r\Bigr)^{1/r},
    \end{aligned}
\end{equation}
where the supremum is taken over all finite length sequences $t_1 \leq t_2 \leq \ldots \leq t_{L}$ in $\mathbb I$. We note that if $\|V^r(F_t:t \in \mathbb I)\|_{L^p(X)} < \infty$ for some $p$, then $V^r(F_t(x):t \in \mathbb I) < \infty$ a.e. $x$ and hence for these $x$, $\lim_{\substack{t \to t_0 \\ t \in \mathbb I}} F_t(x)$ exists for any $t_0.$ The Carleson--Hunt theorem was upgraded to $r$-variational bounds in \cite{osttw}. Let $r'$ be the conjugate exponent of $r$. We have the following theorem.
\begin{theorem}( $r$-variational Carleson--Hunt theorem \cite{osttw} ).\label{osttw}
    Suppose that $r>2$ and $r' <p<\infty$. Then
    \begin{equation}\label{varcar1}
    \begin{aligned}
        \|V^r(S_{N}f:N \in \mathbb N)\|_{L^p(\mathbb T)} \leq C_{p,r}\|f\|_{L^p(\mathbb T)}.
        \end{aligned}
    \end{equation}
    At the endpoint $p=r'$, a restricted weak type result holds, namely, for any $f \in L^{r',1}(\mathbb T)$, the function $V^r(S_{N}f:N\in \mathbb N)$ belongs to $L^{r', \infty}(\mathbb T)$. Also $L^p$ bounds fail for $V^2(S_{N}f: N \in \mathbb N).$
\end{theorem}
 Two interesting applications of the variational bounds of partial Fourier sums (or equivalently, partial Fourier integrals) are given in \cite{osttw}. One is to establish a singular integral variant of the Wiener--Wintner theorem from ergodic theory and the other gives an alternative proof of the Christ--Kiselev nonlinear maximal Hausdorff--Young inequality. For details see \cite{osttw}.

 Let $\{F_t:t \in \mathbb I\}$ be a family of complex-valued measurable functions on a measure space $X$ where $\mathbb I \subseteq \mathbb R.$ For $J \in \mathbb N ,$ we define 
$$\mathcal{G}_{J}(\mathbb I):=\{(t_i:i=1,\ldots,J) \subseteq \mathbb I: t_1 <t_2 < \ldots<t_{J}\};$$
the family of all strictly increasing sequences of length $J$ in $\mathbb I$. For $1 \leq r < \infty$ and a sequence $I=(I_{i}:i=1,\ldots,J) \in \mathcal{G}_{J}(\mathbb I),$ the $r$-oscillation semi-norm is defined by
$$O^{r}_{I,J}(F_{t}(x): t \in \mathbb I):= \Bigl(\sum_{j=1}^{J-1}\sup_{t \in [I_j,I_{j+1})}|F_t(x)-F_{I_j}(x)|^r\Bigr)^{1/r}.$$ 
We have $O^r_{I, J}(F_t: t \in \mathbb I) \leq V^r(F_t : t \in \mathbb I)$ for all $I$ and $J$; so $r$-variational bounds imply uniform $r$-oscillation bounds. Thus uniform oscillation bounds are interesting when variational bounds are not available but they are also useful since they are often much easier to establish than variational bounds. Furthermore, it is easy to see that
$$\|\sup_{t \in \mathbb I}|F_t|\|_{L^p(X)} \leq \sup_{J \in \mathbb N \cup \{\infty\}} \sup_{I \in \mathcal{G}_J(\mathbb I)} \|O^{r}_{I,J}(F_{t}: t \in \mathbb I)\|_{L^p(X)}+\|F_{t_0}\|_{L^p(X)}$$
for any $t_0 \in \mathbb I$ and so establishing $L^p$ bounds for maximal, $r$-oscillation and $r$-variation semi-norms are in increasing order of difficulty. Theorem \ref{osttw} shows that $2$-variational bounds for $\{S_Nf: N \in \mathbb N\}$ \textit{fail}. The next result establishes uniform $2$-oscillation bounds for $\{S_N f: N \in \mathbb N\}$.
\begin{theorem}( $2$-oscillation Carleson-Hunt theorem \cite{MSW} ).\label{sur}
    For $p \in [2, \infty),$ we have the uniform $2$-oscillation bound
    $$\sup_{J \in \mathbb{N} \cup \{\infty\}} \sup_{I \in \mathcal{G}_{J}(\mathbb N)}\|O^{2}_{I,J}(S_{N}f: N \in \mathbb N)\|_{L^p(\mathbb T)} \leq C_p \|f\|_{L^p(\mathbb T)}.$$
    If $\mathbb D$ is the set of dyadic numbers, i.e., numbers of the form $2^n$, then we have
    $$\sup_{J \in \mathbb{N} \cup \{\infty\}} \sup_{I \in \mathcal{G}_{J}(\mathbb D)}\|O^{2}_{I,J}(S_{N}f: N \in \mathbb D)\|_{L^p(\mathbb T)} \leq C_p \|f\|_{L^p(\mathbb T)}, \,\,\,\,\, p \in (1,\infty).$$
\end{theorem}
See \cite{RW} and \cite{MSW} for details. In this paper, we give an extension of Theorems \ref{osttw} and \ref{sur} in higher dimensions.

In higher dimensions, there are a variety of Fourier partial sums to consider. Let $f \in L^1(\mathbb T^d)$ and $\{\widehat{f}(\underline{n})\}_{\underline{n} \in \mathbb Z^d}$ be its Fourier coefficients. Let $\|\cdot\|$ be a norm on $\mathbb R^d$ and consider the one parameter partial sums of the Fourier series of $f$ defined by
\begin{equation}\label{pshd}
\begin{aligned}
    S_{N}f(\underline{x}):=\sum_{\|\underline{n}\| \leq N}\widehat{f}(\underline{n})e^{2 \pi i \underline{n}.\underline{x}}, \,\,\,\, \underline{x} \in \mathbb T^d.
    \end{aligned}
\end{equation}
Interesting cases arise when we consider $\|\cdot\|=\|\cdot\|_{l^q_{d}}$, the standard $l^q_{d}$ norms on $\mathbb R^d$. For example, when the norm is the standard Euclidean norm, that is, $q=2$ (more generally, $1<q<\infty$) the pointwise almost everywhere convergence problem for $L^2(\mathbb T^d)$ functions is a major open problem in harmonic analysis. But we do have positive results when $q=1$ and $q=\infty$ as the balls corresponding to these norms are polygons which do not have any curvature. We have the following theorem due to C. Fefferman.

\begin{theorem}( C. Fefferman \cite{F1} ).\label{fef}
    Let $P$ be a convex polytope in $\mathbb R^d$ such that the origin lies in the interior of $P$. Set $\lambda P=\{\lambda \underline{x}: \underline{x} \in P\}$ for $\lambda > 0$. If $f \in L^p(\mathbb T^d)$, $1<p \leq \infty$ and $\{\widehat{f}(\underline{n})\}_{\underline{n} \in \mathbb Z^d}$ are the Fourier coefficients of $f$, then
    $$f(\underline{x})=\lim_{\lambda \to \infty} \sum_{\underline{n} \in \lambda P}\widehat{f}(\underline{n})e^{2 \pi i \underline{n}.\underline{x}} \,\,\,\, \text{for \ almost \ every \ } \underline{x} \in \mathbb T^d.$$
\end{theorem}
A polytope is a generalisation of polygons in higher dimensions. We shall review the theory of polytopes in the next section. If we define $S_{\lambda}f(\underline{x})=\sum_{\underline{n} \in \lambda P} \widehat{f}(\underline{n})e^{2 \pi i \underline{n}.\underline{x}}$, then Theorem \ref{fef} is a consequence of the maximal function bound
\begin{equation}\label{fefmax}
\begin{aligned}
    \|\sup_{\lambda > 0}|S_{\lambda}f|\|_{L^p(\mathbb T^d)} \leq C \|f\|_{L^p(\mathbb T^d)}.
    \end{aligned}
\end{equation}

The purpose of this paper is to strengthen C. Fefferman's maximal bound (\ref{fefmax}) to a strong $r$-variational bound when $r>2$ and a uniform $2$-oscillation bound.

\begin{theorem}\label{dmt}
    Let $P$ be a convex polytope in $\mathbb R^d$ such that the origin lies in the interior of $P$ and let $\lambda P$ be defined as in Theorem \ref{fef}. For $r>2$ and $r'<p<\infty$, we have
        $$\|V^r(S_{\lambda}f:\lambda > 0)\|_{L^p(\mathbb T^d)} \leq C_{p,r} \|f\|_{L^p(\mathbb T^d)}.$$
    At the endpoint $p=r'$, a restricted weak type estimate holds. For $p \in [2, \infty),$ we have the uniform $2$-oscillation estimate
    $$\sup_{J \in \mathbb N \cup \{\infty\}} \sup_{{I \in \mathcal{G}_{J}(\mathbb R^{+})}}\|O^{2}_{I,J}(S_{\lambda}f: \lambda  > 0)\|_{L^p(\mathbb T^d)} \leq C_{p}\|f\|_{L^p(\mathbb T^d)}.$$ 
    For dyadic indices $\mathbb D=\{2^n:n \in \mathbb Z\}$ and for all $1 < p < \infty,$ we have
    $$\sup_{J \in \mathbb N \cup \{\infty\}} \sup_{{I \in \mathcal{G}_{J}(\mathbb D)}}\|O^{2}_{I,J}(S_{\lambda}f: \lambda \in \mathbb D)\|_{L^p(\mathbb T^d)} \leq C_{p}\|f\|_{L^p(\mathbb T^d)}.$$ 
    \end{theorem}
    
We shall also prove the continuous version of this theorem. Let $f$ be a Schwartz function defined on $\mathbb R^d$ and $\widehat{f}(\underline{\xi})=\int_{\mathbb R^d}f(\underline{x})e^{-2 \pi i \underline{\xi}.\underline{x}} d\underline{x}$ be its Fourier transform. We consider the partial (inverse) Fourier integral
\begin{equation}\label{pip}
\begin{aligned}
    S_{\lambda}f(\underline{x}):=\int_{\underline{\xi} \in \lambda P} \widehat{f}(\underline{\xi}) e^{2 \pi i \underline{x}. \underline{\xi}} d\underline{\xi}.
    \end{aligned}
\end{equation}

\begin{theorem}\label{mt}
    Suppose that $r>2$. Then for $r'<p<\infty$, we have that 
    \begin{equation}\label{varcar2}
    \begin{aligned}
        \|V^r(S_{\lambda}f:\lambda > 0)\|_{L^p(\mathbb R^d)} \leq C_{p,r} \|f\|_{L^p(\mathbb R^d)}.
        \end{aligned}
    \end{equation}
    At the endpoint $p=r'$, a restricted weak type estimate holds. For $p \in [2, \infty),$ we have the uniform $2$-oscillation estimate
    \begin{equation}\label{2osc}
    \begin{aligned}
      \sup_{J \in \mathbb N \cup \{\infty\}} \sup_{{I \in \mathcal{G}_{J}(\mathbb R^{+})}}\|O^{2}_{I,J}(S_{\lambda}f:\lambda > 0 )\|_{L^p(\mathbb R^d)} \leq C_{p}\|f\|_{L^p(\mathbb R^d)}.  
      \end{aligned}
    \end{equation}
    For dyadic indices $\mathbb D=\{2^n: n \in \mathbb Z\}$ and for all $1 < p < \infty,$ we have
    $$\sup_{J \in \mathbb N \cup \{\infty\}} \sup_{{I \in \mathcal{G}_{J}(\mathbb D)}}\|O^{2}_{I,J}(S_{\lambda}f:\lambda \in \mathbb D)\|_{L^p(\mathbb R^d)} \leq C_{p}\|f\|_{L^p(\mathbb R^d)}.$$ 
\end{theorem}
Theorem \ref{mt} implies Theorem \ref{dmt} by a transference argument given in an appendix of \cite{osttw}. See Remark \ref{transference} for further details.

We also consider Fourier partial sums over rectangles. Consider a sequence $\{\mathcal{M}_k\}_{k=1}^{\infty}=\{(M^1_{k},\ldots,M^d_{k})\}_{k=1}^{\infty}$ of $d$-tuples of positive real numbers and define the corresponding axes parallel rectangles $$\mathcal{R}_{k}=\{\underline{x}=(x_1,\ldots,x_d) \in \mathbb R^d: |x_i| \leq M^{i}_k, i=1, \ldots, d\}.$$
For $f \in L^1(\mathbb T^d)$ we consider the rectangular partial sums
$$S_{\mathcal{R}_k}f(\underline{x}): =\sum_{\underline{n} \in \mathcal{R}_k}\widehat{f}(\underline{n})e^{2 \pi i \underline{n}.\underline{x}}, \,\,\,\, \underline{x} \in \mathbb T^d.$$ 
In \cite{F}, C. Fefferman showed that there exists a continuous function $f \in C(\mathbb T^d)$ and a sequence of rectangles $\{\mathcal{R}_k\}$ such that the partial sums $S_{\mathcal{R}_k}f(\underline{x})$ diverge as $k \to \infty$ for \textit{every} $x \in \mathbb T^d$! On the other hand, if the rectangles $\{\mathcal{R}_k\}$ in $\mathbb R^2$ are nested, that is $\mathcal{R}_{k} \subseteq \mathcal{R}_{k+1}$ for all $k$, then we have the following theorem.
\begin{theorem}(Tevzadze \cite{T} and Fefferman \cite{F1} ).\label{tevfef}
   Let $\{\mathcal{R}_k\}_{k \in \mathbb N}$ be a nested sequence of axes parallel rectangles in $\mathbb R^2$. If $f \in L^2(\mathbb T^2),$ then we have $$\|\sup_{k \in \mathbb N}|S_{\mathcal{R}_k}f|\|_{L^{2}(\mathbb T^2)} \leq C\|f\|_{L^2(\mathbb T^2)}.$$ 
   In particular, $\lim_{k \to \infty} S_{\mathcal{R}_k}f(\underline{x})=f(\underline{x})$ a.e. $\underline{x} \in \mathbb T^2$, for all $f \in L^2(\mathbb T^2)$.
\end{theorem}
This was first established by N. Tevzadze in \cite{T}. While his approach is conceptually similar to Fefferman’s technique for proving (\ref{fefmax}) for functions on $\mathbb T^2$, it does not naturally extend to higher dimensions. We shall now show that Fefferman's technique is robust enough to give variation and oscillation bounds in higher dimensions.  We state and prove the continuous version; the discrete version follows from standard transference arguments, which is addressed briefly in Remark \ref{transference}. For $f \in \mathcal{S}(\mathbb R^d)$, consider the partial Fourier integral
$$S_{\mathcal{R}_k}f(\underline{x}):=\int_{\underline{\xi} \in \mathcal{R}_k}\widehat{f}(\underline{\xi})e^{2 \pi i \underline{x}.\underline{\xi}} d\underline{\xi}.$$
\begin{theorem}\label{mt1}
   For a nested sequence of axes parallel rectangles $\{\mathcal{R}_{k}\}$ and for $r>2$ we have 
   \begin{equation}\label{rvarrec}
   \begin{aligned}
     \|V^r(S_{\mathcal{R}_k}f: k \in \mathbb N)\|_{L^2(\mathbb R^d)} \leq C_r \|f\|_{L^2(\mathbb R^d)}  
     \end{aligned}
   \end{equation}
   and
   \begin{equation}\label{2oscrec}
   \begin{aligned}
      \sup_{J \in \mathbb{N} \cup \{\infty\}} \sup_{I \in \mathcal{G}_{J}(\mathbb N)}\|O^{2}_{I,J}(S_{\mathcal{R}_k}f: k \in \mathbb N)\|_{L^2(\mathbb R^d)} \leq C \|f\|_{L^2(\mathbb R^d)}. 
      \end{aligned}
   \end{equation}
\end{theorem}

\subsection*{Acknowledgements} I would like to thank my PhD advisor, Prof. Jim Wright, for insightful discussions and helpful suggestions throughout this project. I am also thankful to the referee for their careful reading of the manuscript and for their insightful comments, which greatly improved the accuracy and clarity of the proofs.

\section{Preliminaries}
\subsection{Polytopes}
We recall the definition and some properties of convex polytopes from \cite{Z}. Although some of the results are well-established in the literature, for instance the triangulation part, we establish its proof here for the sake of completeness.
\begin{definition}
A convex polyhedron $P$ in $\mathbb R^d$ is the intersection of finitely many half-spaces; more precisely, it is the set of all points $\underline{x} \in \mathbb R^d$ that satisfy the system of inequalities
\begin{equation}\label{cp}
\begin{aligned}
  \mathbf{A}\underline{x} \leq \mathbf{b}, 
  \end{aligned}
\end{equation}
where $\mathbf{A}=(a_{ij}) \in \mathcal{M}_{m \times d}(\mathbb R) $ and $ \mathbf{b}=(b_i) \in \mathcal{M}_{m \times 1}(\mathbb R)$.    
\end{definition}
Every row of $\mathbf{A}$ corresponds to a closed half-space $C_i$ whose boundary is a hyperplane, say $H_i$, $i \in \{1,\ldots,m\}$. Hyperplanes $H_i$ are known as bounding hyperplanes and they support $P$ in the sense that $P \cap H_i \neq \emptyset$ and $P$ lies on one side of $H_i$. We shall assume that representation (\ref{cp}) is irredundant; that is, we cannot make the same polyhedron $P$ using fewer inequalities than that given in (\ref{cp}). The representation (\ref{cp}) of $P$ in terms of half-spaces is known as the $\mathcal{H}$-representation of $P$ and we say that $P$ is an $\mathcal{H}$-polyhedron. As $P$ is convex by definition, we speak of polyhedrons without including the word `convex'. It is straight-forward to show that the boundary of $P$ is 
\begin{equation} \label{bdryP}
\begin{aligned}
  \partial P=\cup_{i=1}^{m}(H_i \cap P);  
  \end{aligned}
\end{equation}
and for each $i \in \{1,\ldots,m\}$, $H_i \cap P$ is also a polyhedron, known as a facet of $P$. We have the following notion of a general face of a polyhedron.
\begin{definition}
For a subset $I \subseteq \{1,\ldots,m\}$, the set $$F :=\{\underline{x} \in P: a_i^{T}\underline{x}=b_i, i \in I \}$$ is said to be a face of $P$. Here, for fixed $i$, $a_i^{T}=(a_{ij}) \in \mathcal{M}_{1 \times d}(\mathbb R)$ is the $i^{\text{th}}$ row of $\mathbf{A}$ in (\ref{cp}).  
\end{definition}

 In other words, a face of $P$ is the intersection of $P$ with some bounding hyperplanes. Clearly, every face of $P$ is also a polyhedron. We say that a face $F$ of $P$ has dimension $k$ if the dimension of the smallest affine space containing $F$ is $k$. The (affine) dimension of an affine space $A$ is the dimension of the associated vector space obtained by translating $A$ by a fixed vector in $A$. Notice that if the interior of $P$ is not empty, then a facet is a $(d-1)$ dimensional face as it lies on a hyperplane whose dimension is $(d-1)$. A $0$-dimensional face is known as a vertex of $P$. 
 
 A bounded polyhedron is known as a polytope and polytopes have a simple representation in terms of the vertices which is equivalent to the representation (\ref{cp}).
 We recall that the convex hull of a set $Y \subseteq \mathbb R^d$ is the smallest convex set containing $Y$ and it is denoted by $\text{conv}(Y)$. Let $P$ be a polytope and let $V$ be the set of vertices of $P$, that is, the set of $0$-dimensional faces of $P$. The famous Krein--Milman theorem (see \cite{R}, Chapter 3) says that $V$ is a finite set and
\begin{equation}\label{cpv}
\begin{aligned}
  P=\text{conv} (V).
  \end{aligned}
\end{equation}
The representation of polytope $P$ as the convex hull of its vertices is known as the $\mathcal{V}$-representation of $P$ and we say that $P$ is a $\mathcal{V}$-polytope. We shall also need the following notion of a cone.
\begin{definition}
    A nonempty set of vectors $C \subseteq \mathbb R^d$ is called a cone if for every finite subset of $C$, the set of all corresponding linear combinations with nonnegative coefficients also lies in $C$.
\end{definition}
 We define the conical hull (or positive hull) of a set $Y \subseteq \mathbb R^d$ to be the intersection of all cones in $\mathbb R^d$ that contain $Y$; it is denoted by $\text{cone}(Y)$. One can easily see that
 $$\text{cone}(Y)=\{t_1\underline{y_1}+\ldots+t_n\underline{y_n}:\{\underline{y_1},\ldots,\underline{y_n}\} \subseteq Y, t_i \geq 0 \}.$$  The relation (\ref{cpv}) implies that
\begin{equation}\label{coneofv}
\begin{aligned}
    \text{cone}(P)=\text{cone}(V).
    \end{aligned}
\end{equation}
This follows from the fact that every vector in $\text{cone}(P)$ is a conical combination of vectors in $P$ and the vectors in $P$ are convex combinations of vectors in $V$; we obtain the forward inclusion by unpacking the terms. The reverse inclusion is trivial. This is a crucial observation as it says that the conic hull of $P$ is the same as the conic hull of a finite set.

\subsection{Triangulation of $P$}
Triangulating a polytope with respect to a point is a way of decomposing the polytope into finitely many polytopes, each corresponding to a facet of $P$; this decomposition method allows us to use known one dimensional results for proving Theorem \ref{mt}.
 Let $P$ be a polytope in $\mathbb R^d$ whose $\mathcal{H}$- and $\mathcal{V}$-representation are given by (\ref{cp}) and (\ref{cpv}) respectively. We also assume that the origin is an interior point of $P$. Fix $i \in \{1,\ldots,m\}$. For the facet $F_i:= (P \cap H_i)$ of $P$, let
\begin{equation}\label{Pidefn}
\begin{aligned}
    P_i := \{t \underline{x}: \underline{x} \in F_i, 0 \leq t \leq 1\}=\text{conv}(\{\mathbf{0}\} \cup F_i).
    \end{aligned}
\end{equation} 
We have the following lemma.
\begin{lemma}\label{tri}
Let $P$ be a polytope in $\mathbb R^d$ such that the origin lies in the interior of $P$. For each $i \in \{1,\ldots,m\}$, we construct $P_i$ as in (\ref{Pidefn}). Then we have the following.
\begin{enumerate}
    \item Each $P_i$ is bounded and $P= \bigcup_{i=1}^{m}P_i$;
    \item If $P_i \cap P_j \neq \emptyset$ for some $i \neq j$, then $P_i \cap P_j$ is contained in a hyperplane, hence the $d$-dimensional Lebesgue measure $\mu(P_i \cap P_j)=0.$ 
\end{enumerate}
\end{lemma}

\begin{proof}
\begin{enumerate}
    \item Let $\underline{y} \in P_i$ for some $i$. Then $\underline{y}= t \underline{x}$ for some $\underline{x} \in F_i$ and for some $0 \leq t \leq 1$. As $P$ is convex and $\underline{x}$ and the origin lies in $P$, $\underline{y}=(1-t)\mathbf{0} + t \underline{x} \in P$. This means that $P_i \subseteq P$ and because $P$ is bounded, $P_i$ is also bounded. Let $\underline{y} \in P$. Consider the ray $\overrightarrow{r}(t)=t \underline{y}$, $t > 0$, originating from the origin and passing through $\underline{y}$. As $P$ is closed and bounded, there is some $t_0 \geq 1$ such that $t_0 \underline{y} \in \partial P$. By (\ref{bdryP}), $t_0 \underline{y} \in H_i \cap P= F_i$ for some $i \in \{1,\ldots,m\}$, thus $\underline{y}=\frac{1}{t_0}(t_0 \underline{y}) \in P_i$.
    
    \item We claim that $P_i \cap P_j=\text{conv}(\{\mathbf{0}\} \cup (F_i \cap F_j)).$ Let $\mathbf{0} \neq \underline{x} \in P_i \cap P_j.$ Then there is $t_i, t_j \in (0,1]$, $\underline{u} \in F_i$ and $\underline{v} \in F_j$ such that 
    \begin{equation}\label{uandv}
    \begin{aligned}
        \underline{x}=t_i \underline{u}= t_j \underline{v}.
        \end{aligned}
    \end{equation}
Thus $\underline{u}=t \underline{v}$ where $t=t_j/t_i>0$. Consider the Minkowski functional $\gamma_{P}$ related to $P$ given by $\gamma_P(\underline{y})=\inf \{s>0: \underline{y} \in sP\}.$
    We know that $\gamma_{P}(c\underline{y})=c \gamma_{P}(\underline{y})$ for all $c \geq 0$ and $\gamma_{P}(\underline{y})=1$ if and only if $\underline{y} \in \partial P.$ Therefore, 
    $$1=\gamma_{P}(\underline{u})=\gamma_{P}(t\underline{v})= t \gamma_{P}(\underline{v})=t,$$
    which shows that $t_i=t_j$ and hence $\underline{u}=\underline{v}$. Therefore by (\ref{uandv}), $\underline{x}$ is a convex combination of $\mathbf{0}$ and an element $\underline{u}=\underline{v} \in F_i \cap F_j$; this shows the forward inclusion. The reverse inclusion is immediate. 
    
    Note that $F_i \cap F_j$ is either the empty set or a polytope of dimension at most $(d-2)$, so $P_i \cap P_j =\text{conv}(\{\mathbf{0}\} \cup (F_i \cap F_j))$ is a polytope of dimension at most $(d-1)$ and hence it lies on a hyperplane.
\end{enumerate}
\end{proof}
\begin{remark}\label{decompoly}
\begin{enumerate}
    \item The triangulated components of $P$ are compatible with the dilation structure, that is, for fixed $i =1,\ldots,d$, if $\lambda_1 \leq \lambda_2$, then $\lambda_1 P_i \subseteq \lambda_2 P_i$; for every $\lambda>0$, $\lambda P=\cup_{i=1}^{m} \lambda P_i$. By the previous lemma, the partial Fourier integral of $f$ over $P$ is
    \begin{equation}\label{decomp}
    \begin{aligned}
        S_{\lambda}f(\underline{x})=\int_{\underline{\xi} \in \lambda P}\widehat{f}(\underline{\xi})e^{2 \pi i \underline{x}. \underline{\xi}} d \underline{\xi}= \sum_{i=1}^{m} \int_{\underline{\xi} \in \lambda P_i}\widehat{f}(\underline{\xi})e^{2 \pi i \underline{x}. \underline{\xi}} d \underline{\xi}.
        \end{aligned}
    \end{equation}
    Therefore, it is sufficient to prove Theorem \ref{mt} for a triangulated component of $P$.
    \item For each facet $F_i$, $i \in \{1,\ldots,m\}$, consider the set
 \begin{equation}\label{sector}
 \begin{aligned}
  S_i :=\{\lambda \underline{x}: \underline{x} \in F_i, \lambda \geq 0\}.
  \end{aligned}
 \end{equation}
 The convexity of $F_i$ implies that $S_i=\text{cone}(F_i)$, hence $S_i$ is a cone associated to $F_i$. If $V_i$ is the set of vertices of $F_i$, then by (\ref{coneofv}) we have that $S_i=\text{cone}(V_i)$, hence $S_i$ is a finitely generated cone. By a standard theorem on cones (see for example, \cite{Z}, Theorem 1.3), $S_i$ is the intersection of finitely many half-spaces and therefore, it is an $\mathcal{H}$-polyhedron. We also observe that
 \begin{equation}\label{Pinew}
 \begin{aligned}
     P_i= S_i \cap C_i;
     \end{aligned}
 \end{equation}
 recall that $C_i$ are the closed half-spaces whose intersection is $P$. This means that $P_i$ is a part of the cone $S_i$ obtained by cutting it by the hyperplane $H_i=\{\underline{x} \in \mathbb R^d: a_i^{T}\underline{x} = b_i\}$, hence the name triangular polytope.    
\end{enumerate}
\end{remark}

\subsection{Partitioning of the rectangles}
  As the rectangles $\{\mathcal{R}_k\}_{k}$ are nested, the coordinates $\{\mathcal{M}_k\}_{k}$ corresponding to these rectangles satisfy the coordinatewise order relation $\mathcal{M}_{k} \leq \mathcal{M}_{k+1}$ in $\mathbb R^d$, for all $k \in \mathbb N$. This means that $M^{i}_{k} \leq M^{i}_{k+1}$ for all $i=1,\ldots,d$, and for all $k \in \mathbb N$. We partition these rectangles into finitely many components so that the decomposition is compatible with the decomposition of the rectangles inside it. This decomposition allows us to use one dimensional results to prove Theorem \ref{mt1}. Fix $k \in \mathbb N$. We see that $\mathcal{R}_k$ is the intersection of $2d$ half-planes
$$\mathcal{R}_k=\Bigl(\bigcap_{i=1}^{d}\{\underline{x}: x_i \leq M^i_k\}\Bigr) \cap \Bigl(\bigcap_{i=1}^{d}\{\underline{x}:x_i \geq -M^i_k\}\Bigr).$$
In other words, $\mathcal{R}_k$ is a polytope whose $\mathcal{H}$-representation is given by $\mathbf{A}_{k}\underline{x} \leq \mathbf{b}_{k}$, that is,
 \begin{equation}\label{H-rep}
 \begin{aligned}
\left[ \begin{array}{cccc}
    1 & 0 & \cdots  & 0 \\
    0 & 1 & \cdots  & 0 \\
    \vdots & \vdots & \ddots & \vdots \\
    0 & 0 & \cdots & 1 \\
    \hline
    -1 & 0 & \cdots  & 0 \\
    0 & -1 & \cdots  & 0 \\
    \vdots & \vdots & \ddots & \vdots \\
    0 & 0 & \cdots & -1 \\
\end{array} \right]_{2d \times d}
\begin{bmatrix}
x_1 \\
x_2 \\
\vdots \\
x_d
\end{bmatrix}_{d \times 1}
\leq
\left[ \begin{array}{c}
 M^1_{k}\\
 M^2_{k} \\
\vdots\\
 M^d_{k} \\
 \hline
 M^1_{k} \\
  M^2_{k} \\
\vdots\\
 M^d_{k}
\end{array} \right]_{2d \times 1}
\end{aligned}
\end{equation}
Observe that $\textbf{A}_{k}=\textbf{A}= \begin{bmatrix} I_d ; \, -I_d \end{bmatrix}$ for all $k \in \mathbb N$. It is straight-forward to see that the set of vertices $V_k$ of $\mathcal{R}_{k}$ is $V_{k}=\{(\pm M^1_{k}, \pm M^2_{k}, \ldots, \pm M^d_{k})\}$. Therefore the $\mathcal{V}$-representation of $\mathcal{R}_{k}$ is the convex hull of its vertices 
\begin{equation}\label{V-rep}
\begin{aligned}
\mathcal{R}_{k}=\text{conv} (V_{k}).
\end{aligned}
\end{equation}
Facets of $\mathcal{R}_k$ are
$$F^i_{k}= \{\underline{x} \in \mathcal{R}_k: x_i=M^i_{k}\} \text{ and } \widetilde{F}^i_{k}=\{\underline{x} \in \mathcal{R}_k: x_i=-M^i_{k}\}, \,\,\, i=1,\ldots,d.$$ Every facet itself is a $(d-1)$ dimensional polytope and the set of vertices $V^i_k$ of the facet $F^i_{k}$ is the set of all vertices of $\mathcal{R}_k$ that lie on $F^i_k$. This means that $$V^i_{k}=\{(\pm M^1_{k}, \ldots, M^i_{k}, \ldots, \pm M^d_{k})\}.$$ Similarly, the set of vertices of $\widetilde{F}^i_{k}$ is $\widetilde{V}^i_{k}=\{(\pm M^1_{k}, \ldots, -M^i_{k}, \ldots, \pm M^d_{k})\}.$ Thus we write $F^i_k$ and $\widetilde{F}^i_k$ in its $\mathcal{V}$-representation as $F^i_k=\text{conv}(V^i_k)$ and $\widetilde{F}^i_k=\text{conv}(\widetilde{V}^i_k)$.

We partition each rectangle $\mathcal{R}_k$ into $2d$ components, each component corresponds to a distinct facet. These subregions will be shown to form polyhedral complexes. Recall that a polyhedral complex is a finite collection of convex polytopes satisfying the following conditions:
\begin{enumerate}
    \item The face condition: If a polytope belongs to the complex, then all of its faces must also belong to the complex.
    \item The intersection condition: If any two polytopes in the polyhedral complex intersect, then their intersection must be exactly a single face of both polytopes.
\end{enumerate}
 Later on, we will slightly abuse this definition by calling a polyhedral complex not the collection itself, but the set defined as the union of all its elements.
We begin by triangulating the first rectangle $\mathcal{R}_1$ with respect to the origin. Defining
$$P^i_{1} := \text{conv}(\{\mathbf{0}\} \cup F^i_{1}) \quad \text{and} \quad \widetilde{P}^i_{1} := \text{conv}(\{\mathbf{0}\} \cup \widetilde{F}^i_{1})$$
decomposes the rectangle $\mathcal{R}_1$ into
$$\mathcal{R}_1= \Bigl(\bigcup_{i=1}^{d}P^i_{1}\Bigr) \cup \Bigl(\bigcup_{i=1}^{d}\widetilde{P}^i_{1}\Bigr).$$

To partition $\mathcal{R}_2$ we fix one of its facets, say $F^i_{2}$. Note that the facet $F^i_{2}$ is parallel to the facet $F^i_{1}$ of rectangle $\mathcal{R}_1$. Define 
\begin{equation}\label{conv}
\begin{aligned}
  T^i_{2}:= \text{conv} (F^i_{1} \cup F^i_{2}).
  \end{aligned}
\end{equation}
Because $F^i_1$ and $F^i_2$ are the convex hulls of $V^i_1$ and $V^i_{2}$ respectively, $T^i_{2}$ is the convex hull of $(V^i_{1} \cup V^i_{2})$ and therefore $T^i_2$ is also a convex polytope. We  give an alternative description of $T^i_2$ that will give us more information about the geometry of the set. We claim that 
\begin{equation}\label{sweep}
\begin{aligned}
   T^i_2 = \left\{ \underline{x}: x_i = (1-t)M^i_1 + tM^i_2, \, |x_k| \leq (1-t)M^k_1 + tM^k_2 \, ( k \neq i), \, t \in [0,1] \right\}
   \end{aligned}
\end{equation}
Suppose that $\underline{x} \in \text{conv}(F^i_1 \cup F^i_2)$. Then there is $t \in [0,1]$, $\underline{u} \in F^i_1$ and $\underline{v} \in F^i_2$ such that $\underline{x}=(1-t)\underline{u}+t \underline{v}$. As $\underline{u} \in F^i_1 \subseteq \mathcal{R}_1$, we have $u_i=M^i_1$ and for $k \neq i$, $|u_k| \leq M^k_1$. Similarly, $v_i=M^i_2$ and $|v_k| \leq M^k_2$ for $k \neq i.$ By triangle inequality, $|x_k| \leq (1-t)M^k_1+tM^k_2$. On the other hand, we assume that $\underline{x}$ satisfies 
$$x_i=(1-t)M^i_1+tM^i_2,\,\, |x_k| \leq (1-t)M^k_1+tM^k_2,\, ( k \neq i), \text{ for some } t \in [0,1].$$ For $k \neq i$, define the ratio
$$\theta_k=\frac{x_k}{(1-t)M^k_1+tM^k_2}.$$
Clearly $|\theta_k| \leq 1.$ Define $\underline{u}=(u_1,\ldots,u_d)$ and $\underline{v}=(v_1,\ldots,v_d)$ as
$$u_i=M^i_1,\,\, v_i=M^i_2, \,\, u_k=\theta_k M^k_1 \text{  and  }  v_k=\theta_k M^k_2, \, k \neq i.$$
Clearly, $\underline{u} \in F^i_1$ and $\underline{v} \in F^i_2$ and $\underline{x}=(1-t)\underline{u}+t\underline{v}$. This proves our claim.

Now we define
\begin{equation}\label{P2}
\begin{aligned}
  P^i_2:= P^i_1 \cup T^i_2.  
  \end{aligned}
\end{equation}
Observe that $P^i_1$ and $T^i_2$ lie on the either side of the hyperplane $\{\underline{x}: x_i=M^i_1\}$ and their intersection is exactly the facet $F^i_1$, hence $P^i_2$ is a polyhedral complex. We define $\widetilde{P}^i_2$ similarly using $\widetilde{P}^i_1$ and facets $\widetilde{F}^i_1$ and $\widetilde{F}^i_2$. We have the following lemma.

\begin{lemma}\label{parrec}
    The rectangle $\mathcal{R}_2$ satisfies:
    \begin{enumerate}
        \item $\mathcal{R}_2= \Bigl(\bigcup_{i=1}^{d}P^i_{2}\Bigr) \cup \Bigl(\bigcup_{i=1}^{d}\widetilde{P}^i_{2}\Bigr),$ where $P^i_2$ and $\widetilde{P}^i_2$ are defined as in (\ref{P2}) ; 
        \item If $i \neq j$, the intersection $P^i_2 \cap P^j_2$ is contained in a finite union of hyperplanes and has $d$-dimensional Lebesgue measure zero. The same conclusion holds for the intersections $\widetilde{P}^i_2 \cap \widetilde{P}^j_2$, $\widetilde{P}^i_2 \cap P^j_2$, and $P^i_2 \cap \widetilde{P}^j_2$.
    \end{enumerate}
\end{lemma}

\begin{proof}
    \begin{enumerate}
        \item Let $\underline{x} \in \mathcal{R}_2$. If $\underline{x} \in \mathcal{R}_1$, then $\underline{x} \in P^j_1$ (or $\widetilde{P}^j_1$) for some $j$ and we are done. Let $\underline{x} \in \mathcal{R}_2 \setminus \mathcal{R}_1.$ For each $i \in \{1, \ldots, d\}$ there is $t_i \in (-\infty, 1]$ such that $|x_i|=(1-t_i)M^i_1+t_i M^i_2$. Let $t_j=\max_{i} t_i.$ Because $\underline{x} \in \mathcal{R}_2 \setminus \mathcal{R}_1$, $t_j \geq 0$. For any $i \neq j$, we have $t_i \leq t_j$ and hence $$|x_i|=(1-t_i)M^i_1 +t_i M^i_2 \leq (1-t_j)M^i_1+t_j M^i_2.$$
        Thus $\underline{x} \in T^j_2$ or $ \underline{x} \in \widetilde{T}^j_2$ depending on weather $x_j=(1-t_j)M^j_1+t_jM^j_2$ or $x_j=(1-t_j)(-M^j_1)+t_j(-M^j_2).$
        This proves the forward inclusion. The reverse inclusion is immediate as the sets $P^i_2$ and $\widetilde{P}^i_2$ are contained in $\mathcal{R}_2$ by their construction.
        \item Without loss of generality we consider sets $P^i_2$ and $P^j_2$, $i \neq j$, and observe that
        $$P^i_2 \cap P^j_2=(P^i_1 \cap P^j_1) \cup (P^i_1 \cap T^j_2) \cup (P^j_1 \cap T^i_2) \cup (T^i_2 \cap T^j_2).$$
        We have already proved in Lemma \ref{tri} that the intersection of any two triangulated components of a polytope lie in a hyperplane, so the Lebesgue measure of the intersection $\mu(P^i_1 \cap P^j_1)=0.$  For the second term, we observe that $T^j_2$ lies entirely in the half-plane $\{\underline{x}:x_j \geq M^j_1\}$ and $P^i_1$ lies entirely in the half-plane $\{\underline{x}: x_j \leq M^j_1\}.$ Thus their intersection is a measurable set lying in the hyperplane $\{\underline{x}: x_j=M^j_1\}$ and hence $\mu(P^i_1 \cap T^j_2)=0$. Similarly, $\mu(P^j_1 \cap T^i_2)=0.$ 

        It remains to show that $\mu(T^i_2 \cap T^j_2)=0.$ Consider $F^{ij}_1=F^i_{1} \cap F^j_{1}$ and $F^{ij}_2=F^{i}_2 \cap F^{j}_2$; they are convex polytopes each of dimension $(d-2).$ Define $\Sigma^{ij}=\text{conv}(F^{ij}_1 \cup F^{ij}_2).$ Notice that $\Sigma^{ij}$ lies on a hyperplane. We claim that 
        \begin{equation}
        \begin{aligned}
         T^i_2 \cap T^j_2=\Sigma^{ij}.   
         \end{aligned}
        \end{equation}
         Let $\underline{x} \in \Sigma^{ij}.$ Then $\underline{x}$ is a convex combination of elements in $F^{ij}_1 \cup F^{ij}_2.$ As $F^{ij}_1 \subseteq F^i_1$ and $F^{ij}_2 \subseteq F^{i}_2$, $\underline{x}$ is a convex combination of elements in $F^i_1 \cup F^i_2$ and hence $\underline{x} \in T^i_2$. Similarly, $\underline{x} \in T^j_2$ and therefore $\underline{x} \in T^i_2 \cap T^j_2.$ 
        
        Let $\underline{x} \in T^i_2 \cap T^j_2.$ This means that there is $t_i, t_j \in [0,1], \underline{u} \in F^{i}_1, \underline{v} \in F^i_{2},  \underline{w} \in F^j_1$ and $\underline{z} \in F^j_2$ such that 
        \begin{equation}
        \begin{aligned}
            \underline{x}=(1-t_i)\underline{u}+t_i\underline{v} \text{  and  } \underline{x}=(1-t_j)\underline{w}+t_j\underline{z}.
            \end{aligned}
        \end{equation}
        Thus the $j^{\text{th}}$ coordinate $x_j$ of $\underline{x}$ satisfies $x_j=(1-t_i)u_j+t_iv_i \leq (1-t_i)M^j_1+t_i M^j_2$; this gives an upper bound for $x_j.$ We also note that $x_j=(1-t_j)w_j+t_jz_j=(1-t_j)M^j_1+t_j M^j_2$ because $\underline{w} \in F^j_1$ and $\underline{z} \in F^j_2$. If $t_i < t_j$ then we have
        $$x_j \leq (1-t_i)M^j_1+t_i M^j_2 < (1-t_j)M^j_1+t_j M^j_2=x_j,$$ a contradiction. Similarly, we can show that $t_j<t_i$ is not possible. Thus $t_i=t_j=t.$ This allows us to conclude that $\underline{u}, \underline{w} \in F^{ij}_1$ and $\underline{v}, \underline{z} \in F^{ij}_2$ so that $\underline{x}=(1-t)(\frac{1}{2}\underline{u}+\frac{1}{2}\underline{w})+t(\frac{1}{2}\underline{v}+\frac{1}{2}\underline{z})$ with $\frac{1}{2}\underline{u}+\frac{1}{2}\underline{w} \in F^{ij}_1$ and $\frac{1}{2}\underline{v}+\frac{1}{2}\underline{z} \in F^{ij}_2.$  Hence $\underline{x} \in \Sigma^{ij}.$ But $\Sigma^{ij}$  sits inside a hyperplane, so $\mu(\Sigma^{ij})=\mu(T^i_2 \cap T^j_2)=0.$ We can make appropriate changes in the proof to show that $\mu(\widetilde{P}^i_2 \cap P^j_2)= \mu(P^i_2 \cap \widetilde{P}^j_2)= \mu(\widetilde{P}^i_2 \cap \widetilde{P}^j_2)=0$. This proves the lemma.
    \end{enumerate}
\end{proof}

\begin{remark}\label{remrec}
For every $k \geq 3,$ we construct polytopes $T^i_{k}$ and $\widetilde{T}^i_{k}$ as in (\ref{conv}) and define the polyhedral complexes $P^i_k$ and $\widetilde{P}^i_k$ as in (\ref{P2}) that decompose $\mathcal{R}_k$ as
 \begin{equation}
 \begin{aligned}
    \mathcal{R}_{k}=\Bigl(\bigcup_{i=1}^{d}P^i_{k}\Bigr) \cup \Bigl(\bigcup_{i=1}^{d}\widetilde{P}^i_{k}\Bigr) 
    \end{aligned}
 \end{equation} 
 such that for every $i \neq j$, $\mu(P^i_k \cap P^j_k)= 0$ and the same conclusion holds for the intersections $\widetilde{P}^i_k \cap \widetilde{P}^j_k$, $\widetilde{P}^i_k \cap P^j_k$, and $P^i_k \cap \widetilde{P}^j_k$.
 Moreover, arbitrary union of $P^i_k \text{ and/or } \widetilde{P}^i_k$ forms a polyhedral complex. In particular, above decomposition is a way of looking $\mathcal{R}_k$ as a polyhedral complex.
 One of the important features of this partition is that for every $i=1,\ldots,d$, we have the inclusion relation 
 \begin{equation}\label{nested}
 \begin{aligned}
   P^i_{1} \subseteq P^i_{2} \subseteq \ldots \text{ and } \widetilde{P}^i_{1} \subseteq \widetilde{P}^i_{2} \subseteq \ldots.  
   \end{aligned}
 \end{equation}
For $i=1,\ldots,d$, define the infinite polyhedral complexes 
\begin{equation}
\begin{aligned}
 \mathcal{P}^i:= \bigcup_{k=1}^{\infty} P^i_{k}\,\, \text{ and } \,\, \widetilde{\mathcal{P}}^i:= \bigcup_{k=1}^{\infty} \widetilde{P}^i_{k}   
 \end{aligned}
\end{equation}
By the extension of Lemma \ref{parrec} from $\mathcal{R}_2$ to $\mathcal{R}_k$, we have
\begin{equation}\label{decomprec}
\begin{aligned}
  S_{\mathcal{R}_k}f(\underline{x})= \sum_{i=1}^{d} S_{P^i_k}f(\underline{x}) + \sum_{i=1}^{d} S_{\widetilde{P}^i_k}f(\underline{x}).  
  \end{aligned}
\end{equation}
 Also because (\ref{nested}) holds, it is sufficient to prove Theorem \ref{mt1} for one such sequence $\{P^i_k\}_{k \in \mathbb N}$ of polyhedral complexes.
\end{remark}

\section{Proof of the Main theorems}
\begin{proof}[Proof of Theorem \ref{mt}]
We shall prove the restricted weak type estimate
\begin{equation}\label{wte}
\begin{aligned}
    \|V^r(S_{\lambda}f:\lambda > 0)\|_{L^{p,\infty}(\mathbb R^d)} \leq C_{p,r}\|f\|_{L^{p,1}(\mathbb R^d)},
    \end{aligned}
\end{equation}
for $r' \leq p < \infty$ which implies the strong type bounds (\ref{varcar2}) for $r' <p<\infty$. By Remark \ref{decompoly} it suffices to prove the theorem for a triangular polytope $P$.

Let $F$ be the facet of $P$ opposite to the origin lying on the hyperplane $H$ corresponding to the half-space $C$ and let $S$ be the cone associated to $F$. Without loss of generality, we may replace $P$ by a suitable dilate, still denoted by $P$ such that that $\text{dist}(\textbf{0}, H)=1$. Then there is a unique $\underline{n} \in H$ such that $\underline{n}$ is a unit vector. Let $R$ be the rotation linear transformation such that the image of $\underline{n}$ under $R$ is $e_1=(1,0,\ldots,0)$, the standard unit vector along the $x_1$-axis. Thus the image $R(P)=\{t \underline{x}: \underline{x} \in R(F), 0 \leq t \leq 1\}$ is the rotated polytope so that its priciple facet $R(F)$ lies in the hyperplane $\{\underline{x}: x_1=1\}$. We also observe from the cone-half-plane representation (\ref{Pinew}) of a polytope that the image of $P$ under $R$ is
$$R(P)= R(S) \cap R(C)=R(S) \cap \{\underline{x} \in \mathbb R^d: x_1 \leq 1\}=\{\underline{x} \in R(S):x_1 \leq 1\}.$$
Because $R$ is a unitary transformation, it is sufficient to prove the theorem when $P=\{\underline{x} \in S:x_1 \leq 1\}$. Let $f$ be a Schwartz function on $\mathbb R^d$. If we write $\underline{x}=(x_1,x_2,\ldots,x_d)=(x_1,x')$ and $\underline{\xi}=(\xi_1,\xi_2,\ldots,\xi_d)=(\xi_1,\xi')$, where $x',\xi' \in \mathbb R^{d-1}$, then 
    \begin{align*}
        S_{\lambda}f(\underline{x}) &= \int_{\lambda P}\widehat{f}(\underline{\xi})e^{2 \pi i \underline{x}.\underline{\xi}}d\underline{\xi} \\
        &= \int_{\xi_1=0}^{\lambda}\Bigl(\int_{\xi':\underline{\xi} \in S}\widehat{f}(\underline{\xi})e^{2 \pi i x'.\xi'}d \xi'\Bigr)e^{2 \pi i x_1 \xi_1}d \xi_1 \\
        &= \int_{\xi_1=0}^{\lambda} g_{x'}(\xi_1)e^{2 \pi i x_1 \xi_1}d \xi_1,
    \end{align*}
    where $g_{x'}(\xi_1)=\int_{\xi':\underline{\xi} \in S}\widehat{f}(\underline{\xi})e^{2 \pi i x'.\xi'}d \xi'$. If we fix $x' \in \mathbb R^{d-1}$, then we can define a new partial Fourier integral operator
    $$T_{\lambda}g_{x'}(x_1)=\int_{\xi_1=0}^{\lambda}g_{x'}(\xi_1)e^{2 \pi i x_1 \xi_1} d\xi_1,$$
    on functions $g_{x'}$ that rapidly decay as $\xi_1 \to \pm \infty$ and which are obtained from $f \in \mathcal{S}(\mathbb R^d)$ by freezing the last $(d-1)$ variables; the two operators are related by the formula
    \begin{equation}\label{freeze}
    \begin{aligned}
        S_{\lambda}f(\underline{x})=S_{\lambda}f(x_1,x')=T_{\lambda}g_{x'}(x_1).
        \end{aligned}
    \end{equation}
    We see that $T_{\lambda}g_{x'}$ is obtained from $S_{\lambda}f$ by freezing the last $(d-1)$ variables. We apply Theorem \ref{osttw} to obtain
    \begin{equation}\label{apposttw}
    \begin{aligned}
        \|V^r(T_{\lambda}g_{x'}:\lambda > 0)\|_{L^{p,\infty}(\mathbb R)} \leq C_{p,r}\|\check{g}_{x'}\|_{L^{p,1}(\mathbb R)},
        \end{aligned}
    \end{equation}
    where $\check{g}_{x'}$ is the inverse Fourier transform of $g_{x'}$ and $r' \leq p <\infty$. As (\ref{freeze}) holds for every $\lambda > 0$, we have that $V_{r}(S_{\lambda}f(\underline{x}):\lambda> 0)=V_r(T_{\lambda}g_{x'}(x_1):\lambda > 0)$, so the distribution functions satisfy
    $$d_{V_r(S_{\lambda}f:\lambda > 0)}(\alpha)=\int_{\mathbb R^{d-1}}d_{V_r(T_{\lambda}g_{x'}:\lambda > 0)}(\alpha) dx', \,\,\, \alpha \geq 0,$$ where the distribution function of a measurable function $h: X \to \mathbb C$ is given by $d_h(\alpha)=\mu(\{x \in X: |h(x)| \geq \alpha \})$. Therefore, the $L^{p,\infty}$ norms satisfy the inequality
    \begin{equation}\label{weaknorms}
    \begin{aligned}
      \|V^r(S_{\lambda}f:\lambda > 0)\|_{L^{p,\infty}(\mathbb R^d)} \leq \Bigl(\int_{\mathbb R^{d-1}}\|V^r(T_{\lambda}g_{x'}:\lambda>0)\|^p_{L^{p,\infty}(\mathbb R)} dx'\Bigr)^{1/p}.  
      \end{aligned}
    \end{equation}
    Combining the inequalities (\ref{apposttw}) and (\ref{weaknorms}), we have
\begin{align*}
    \|V^r(S_{\lambda}f:\lambda > 0)\|_{L^{p,\infty}(\mathbb R^d)} &\leq C_{p,r}\Bigl(\int_{\mathbb R^{d-1}}\|\check{g}_{x'}\|^p_{L^{p,1}(\mathbb R)} dx' \Bigr)^{1/p} \\
    &= pC_{p,r}\Bigl(\int_{\mathbb R^{d-1}} \Bigl( \int_{0}^{\infty} d_{\check{g}_{x'}}(s)^{1/p} ds\Bigr)^pdx'\Bigr)^{1/p} \\
    &\leq pC_{p,r}\int_{0}^{\infty}\Bigl( \int_{\mathbb R^{d-1}} d_{\check{g}_{x'}}(s) dx'\Bigr)^{1/p} ds.
\end{align*}
    Here, we have expressed the Lorentz norm $\|\check{g}_{x'}\|_{L^{p,1}(\mathbb R)}=p\int_{0}^{\infty}d_{\check{g}_{x'}}(s)^{1/p}ds$ in terms of its distribution function and the last inequality is due to Minkowski integral inequality.
    We observe that
    \begin{align*}
        \check{g}_{x'}(x_1) &= \int_{\mathbb R}g_{x'}(\xi_1)e^{2 \pi i x_1 \xi_1} d\xi_1 \\
        &= \int_{\xi_1 \in \mathbb R} \Bigl( \int_{\xi':\underline{\xi} \in S}\widehat{f}(\underline{\xi})e^{2 \pi i x'.\xi'} d\xi'\Bigr)e^{2 \pi i x_1 \xi_1}d\xi_1 \\
        &= \int_{S}\widehat{f}(\underline{\xi})e^{2 \pi i \underline{x}. \underline{\xi}} d\underline{\xi} \\
        &= \int_{\mathbb R^d}\mathbbm{1}_{S}(\underline{\xi})\widehat{f}(\underline{\xi})e^{2 \pi i \underline{x}.\underline{\xi}} d\underline{\xi}.
    \end{align*}
    Let $Tf(x)=\left(\mathbbm{1}_{S} \, \widehat{f}\right)^{\widecheck{}}(x)$ so that $Tf(x)=\check{g}_{x'}(x_1)$. The distribution functions satisfy $d_{Tf}(s)=\int_{\mathbb R^{d-1}} d_{\check{g}_{x'}}(s)dx'$, which means that
    \begin{align*}
       \|V^r(S_{\lambda}f:\lambda > 0)\|_{L^{p,\infty}(\mathbb R^d)} &\leq p C_{p,r}\int_{0}^{\infty}\Bigl( \int_{\mathbb R^{d-1}} d_{\check{g}_{x'}}(s) dx'\Bigr)^{1/p} ds \\
       &\leq p C_{p,r} \int_{0}^{\infty} d_{Tf}(s)^{1/p} ds = C_{p,r} \|Tf\|_{L^{p,1}(\mathbb R^d)}.
    \end{align*}
    The operator $T$ is a Fourier multiplier of the cone $S$. As observed earlier, $S$ is an $\mathcal{H}$-polyhedron, hence it is the intersection of finitely many half-spaces. Thus, the characteristic function of $S$ is the product of characteristic functions of finitely many half-spaces. As the characteristic function of a half-space is an $L^p$ multiplier for $1 < p< \infty$, and finite product of $L^p$ multipliers is also an $L^p$ multiplier, we have $\|Tf\|_{L^p(\mathbb R^d)} \leq C_{p,S}\|f\|_{L^p(\mathbb R^d)}$, where $C_{p,S}$ is a constant that depends on $p$ and the number of half-spaces $S$ is made up of. Therefore, by the off-diagonal Marcinkiewicz interpolation theorem (see, for instance, \cite{G}, p. 61), we obtain $\|Tf\|_{L^{p,1}(\mathbb R^d)} \leq C'_{p,S}\|f\|_{L^{p,1}(\mathbb R^d)}$ with some new constant $C'_{p,S}$, where we interpolate between carefully chosen end-points that depend only on $p$. This proves the variational bound (\ref{varcar2}). To prove the uniform $2$-oscillation bound (\ref{2osc}), we apply Theorem \ref{sur} to the operators $\{T_{\lambda}: \lambda > 0\}$. This completes the proof of Theorem \ref{mt}.
    \end{proof}
    
\begin{proof}[Proof of Theorem \ref{mt1}]
 Take $f \in \mathcal{S}(\mathbb R^d)$ and write $\underline{x}=(x_1,x')$ and $\underline{\xi}=(\xi_1,\xi')$ as done in the theorem for polytopes. We recall from (\ref{decomprec}) in Remark \ref{remrec} that it is sufficient to prove the theorem for one such sequence of polyhedral complexes; we consider the nested sequence $\{P^1_k\}_{k \in \mathbb N}$ where $\bigcup_{k}P^1_k=\mathcal{P}^1$. Write the partial (inverse) Fourier integrals as
\begin{align*}
    S_{P^1_{k}}f(\underline{x}) &= \int_{\underline{\xi} \in P^1_{k}} \widehat{f}(\underline{\xi}) e^{2 \pi i \underline{x}. \underline{\xi}} d \underline{\xi} \\
    &= \int_{\xi_1=0}^{M^1_{k}} \Bigl( \int_{\xi': \underline{\xi} \in \mathcal{P}^1} \widehat{f}(\underline{\xi}) e^{2 \pi i x'.\xi'} d \xi' \Bigr) e^{2 \pi i x_1 \xi_1} d \xi_1 \\
    &= \int_{\xi_1=0}^{M^1_{k}} g_{x'}(\xi_1) e^{2 \pi i x_1 \xi_1} d \xi_1,
\end{align*}
where $g_{x'}(\xi_1)=\int_{\xi': \underline{\xi} \in \mathcal{P}^1} \widehat{f}(\underline{\xi}) e^{2 \pi i x'.\xi'} d \xi'.$
Note that for fixed $x'$, $g_{x'}$ rapidly decays as $\xi_1 \to \pm \infty$.
For fixed $x' \in \mathbb R^{d-1}$, we define partial Fourier integral operator
$$T_{P^1_k}g_{x'}(x_1):= \int_{\xi_1=0}^{M^1_{k}}g_{x'}(\xi_1)e^{2 \pi i x_1 \xi_1} d \xi_1.$$ Note that
$$S_{P^1_{k}}f(\underline{x})=T_{P^1_{k}}g_{x'}(x_1).$$
Applying Theorem \ref{osttw} to $\{T_{P^1_{k}}g_{x'}: k \in \mathbb N \}$, for $r>2$, we have
\begin{equation}\label{2}
\begin{aligned}
  \|V^r(T_{P^1_{k}}g_{x'}: k \in \mathbb N)\|_{L^2(\mathbb R)} \leq C_r\|\widecheck{g}_{x'}\|_{L^2(\mathbb R)},  
  \end{aligned}
\end{equation}
where $\widecheck{g}_{x'}$ is the inverse Fourier transform of $g_{x'}.$
Observe that $$\widecheck{g}_{x'}(x_1)=\int_{\mathbb R^d}\mathbbm{1}_{\mathcal{P}^1}(\underline{\xi})\widehat{f}(\underline{\xi})e^{2 \pi i \underline{x}.\underline{\xi}} d \underline{\xi}= \left(\mathbbm{1}_{\mathcal{P}^1}\widehat{f}\right)^{\widecheck{}}(\underline{x}),$$ where recall that $\mathcal{P}^1$ is a polyhedral complex and the characteristic function on $\mathcal{P}^1$ is an $L^2$ multiplier. This proves (\ref{rvarrec}). Similarly, we use Theorem \ref{sur} to obtain uniform $2$-oscillation bound (\ref{2oscrec}). This completes the proof of Theorem \ref{mt1}.
\end{proof}

\begin{remark}\label{transference}
    The $r$-variation and $2$-oscillation bounds in the discrete setting are obtained from the estimates in Theorems \ref{mt} and \ref{mt1} using standard transference techniques explained in Appendix A of \cite{osttw}. Suppose that $T$ is a bounded linear operator corresponding to the multiplier $m$ on $\mathbb R^d$, then de Leeuw's transference result \cite{KdL} says that if $m$ is continuous on the lattice points $\mathbb Z^d$, then we can pass bounds of $T$ to its periodised counterpart $\widetilde{T}$. In the context of maximal, variational or oscillation bounds, the framework involves a family of operators defined by a corresponding family of multipliers. To transfer maximal function bounds (see, for example, \cite{G}, p. 281), variational or oscillation bounds, we look at these operators as vector-valued operators. For example, the operator in (\ref{wte}) can be interpreted as a map from $L^{p,1}(\mathbb R^d)$ to $L^{p,\infty}(\mathbb R^d; \ell^{\infty}(\ell^r))$. The vector-valued extension of the transference principle requires that every individual multiplier in this family is continuous on $\mathbb Z^d$.
    
    \begin{enumerate}
        \item \textit{Polytope and its dilates}: As the family of polytopes $\{\lambda P: \lambda>0\}$ is uncountable, if we remove all the dilates $\lambda P$ so that $\partial(\lambda P) \cap \mathbb Z^d \neq \emptyset$, we are still left with a dense set $\Lambda \subseteq \mathbb R^{+}$ such that the boundary of $\lambda P, \lambda \in\Lambda$, does not intersect the lattice points $\mathbb Z^d$. For every $\lambda \in \Lambda$, the multiplier corresponding to $S_{\lambda P}$, that is the characteristic function of $\lambda P$, is continuous on $\mathbb Z^d.$  By the argument in Appendix A of \cite{osttw}, we obtain the $r$-variational bounds for the periodised operator $\|V^r(\widetilde{S}_{\lambda P}f: \lambda \in \Lambda)\|_{L^{p,\infty}(\mathbb T^d)} \leq C_{p,r} \|f\|_{L^{p,1}(\mathbb T^d)}$. Finally, we observe that $V^r(\widetilde{S}_{\lambda P}f:\lambda \in \Lambda)=V^r(\widetilde{S}_{\lambda P}f:\lambda>0)$ for trigonometric polynomials $f$ on $\mathbb T^d$. This gives us the bounds for all the dilates and for trigonometric polynomials. The density of trigonometric polynomials in $L^{p,1}(\mathbb T^d)$ allows us to establish the required bounds on $\mathbb T^d.$
        \item \textit{Sequence of nested rectangles}: For each $k \in \mathbb{N}$, we choose a scaling vector $\underline{\lambda}_k = (\lambda_{1,k}, \dots, \lambda_{d,k}) \in [0,2]^d$ sufficiently close to $\mathbf{1}=(1,\ldots,1)$ and the sequence $\mathcal{R}_{k}'=\{(\lambda_{1,k} x_1,\ldots,\lambda_{d,k} x_d): (x_1,\ldots, x_d) \in \mathcal{R}_k\}$ of nested rectangles so that the characteristic function of $\mathcal{R}_{k}'$ is continuous on $\mathbb Z^d$ and $\mathbbm{1}_{\mathcal{R}_k'}(m) = \mathbbm{1}_{\mathcal{R}_k}(m)$ for all $m \in \mathbb{Z}^d$. An important feature of Theorem \ref{mt1} is that the $r$-variational bound (\ref{rvarrec}) is uniform with respect to \textit{any} sequence of nested rectangles; therefore, the rectangles $\{\mathcal{R}_{k}'\}_{k}$ satisfy the same bound (\ref{rvarrec}). The transference argument in \cite{osttw} allow us to transfer these bounds to the periodised operators $\{\widetilde{S}_{\mathcal{R}_{k}'}: k \in \mathbb N\}$. Finally,  $\mathbbm{1}_{\mathcal{R}_k'}(m) = \mathbbm{1}_{\mathcal{R}_k}(m)$ implies that the partial sum operators are identical on the torus $\mathbb T^d$, that is, $\widetilde{S}_{\mathcal{R}_{k}'} = \widetilde{S}_{\mathcal{R}_{k}}$. Consequently, the uniform $r$-variation and $2$-oscillation bounds hold for $\{\widetilde{S}_{\mathcal{R}_k}: k \in \mathbb N\}$ as well without any change in the bounding constant.
        \end{enumerate}
    
\end{remark}

\end{document}